\documentclass[a4, 11pt]{amsart}
\usepackage[mathscr]{eucal}
\usepackage{amssymb}
\usepackage{latexsym}
\usepackage{amsthm}
\theoremstyle{plain}
\begin{document}
\renewcommand{\theequation}{\arabic{section}.\arabic{equation}}
\newcommand{\be}{\begin{eqnarray}}
\newcommand{\en}{\end{eqnarray}}
\newcommand{\no}{\nonumber}
\newcommand{\la}{\lambda}
\newcommand{\laa}{\Lambda}
\newcommand{\ino}{\int_\Omega}
\newcommand{\tg}{\triangle}
\newcommand{\Li}{\Lambda_i}
\newcommand{\Lj}{\Lambda_j}
\newcommand{\Lk}{\Lambda_{k+1}}
\newcommand{\li}{\lambda_i}
\newcommand{\lj}{\lambda_j}
\newcommand{\lk}{\lambda_{k+1}}
\newcommand{\p}{\partial}
\newcommand{\n}{\nabla}
\newcommand{\nj}{\nabla_j}
\newcommand{\nk}{\nabla_k}
\newcommand{\mc}{\mcite}
\newcommand{\g}{\gamma}
\newcommand{\e}{\varepsilon}
\newcommand{\s}{\sup}
\newcommand{\ep}{\epsilon}
\newcommand{\de}{\delta}
\newcommand{\D}{\Delta}
\newcommand{\pl}{\parallel}
\newcommand{\ov}{\overline}
\newcommand{\bet}{\beta}
\newcommand{\al}{\alpha}
\newcommand{\fr}{\frac}
\newcommand{\pa}{\partial}
\newcommand{\we}{\wedge}
\newcommand{\om}{\Omega}
\newcommand{\na}{\nabla}
\newcommand{\lan}{\langle}
\newcommand{\ra}{\rangle}
\newcommand{\fa}{\sum_{\alpha=1}^n}
\newcommand{\vs}{\vskip0.3cm}
\newcommand{\su}{\sum_{i=1}^k(\Lambda_{k+1}-\Lambda_i)^2}
\newcommand{\suu}{\sum_{i=1}^k(\Lambda_{k+1}-\Lambda_i)}
\renewcommand{\thefootnote}{}
\newcommand{\ri}{\rightarrow}
\title[ The buckling problem of arbitrary  order] {Eigenvalues of  the buckling problem of arbitrary  order on bounded domains of $\mathbb{M}\times\mathbb{R}$}
\footnotetext{2000 {\it Mathematics Subject Classification }: 35P15, 53C20, 53C42, 58G25 \\
 Key words and phrases: Universal inequality for  eigenvalues, the buckling problem of arbitrary order,  ${\mathbb M}\times{\mathbb R}$.}
 \footnotetext{This work was partially supported by CNPq and PROCAD/CAPES.}
\author[Q. Wang and C. Xia]
{ Qiaoling Wang and  Changyu Xia}
\address{Qiaoling Wang \\ \newline \indent
Departamento de Matem\'{a}tica, Universidade de Bras\'{\i}lia, 70910-900-Bras\'{\i}lia-DF, Brazil, wang@mat.unb.br}

\address{Changyu Xia \\  \newline \indent   Departamento de Matem\'atica, Universidade de Bras\'{\i}lia, 70910-900 Bras\'{\i}lia-DF, Brazil.  xia@pq.cnpq.br}
\date{}
\maketitle
\begin{abstract} \noindent
We obtain universal inequalities for  eigenvalues of the buckling problem of arbitrary
order on bounded domains in $\mathbb{M}\times\mathbb{R}$.

\end{abstract}

\section{Introduction}
Let $\om$ be a bounded domain with smooth boundary in an  $n(\geq 2)$-dimensional Riemannian manifold $\mathbb{N}$  and denote by $\Delta $ the Laplace operator acting on functions on $\mathbb{N}$. Let $\nu$ be the outward unit normal vector field of $\pa \om$ and
 let us consider the following  eigenvalue problem :
\begin{eqnarray}\label{eq1}
\left\{\begin{array}{l}(-\Delta)^lu =-\Lambda \D u,~~\mbox{in}
~~\Omega,\\ ~~~~u=\frac{\partial u}{\partial
\nu}=\cdots=\frac{\partial^{l-1} u}{\partial
\nu^{n-1}}=0,~~\mbox{on}~~\partial \Omega,
\end{array}
\right.
\end{eqnarray}
where $l$ is an integer no less than 2. This problem is called the buckling problem of order $l$ which has interpretations in physics, that is, it  describes the critical buckling load of a clamped plate subjected to a uniform
compressive force around its boundary.
Let
\be
\no & & 0<\laa_1\leq\laa_2\leq\laa_3\leq\cdots
\en
denote the successive eigenvalues for (\ref{eq1}).
Here each eigenvalue is repeated according to its
multiplicity. An
important theme of geometric analysis is to estimate these (and other) eigenvalues. When $l=2$ and $\om$ is a bounded domain in an
$n$-dimensional Euclidean space $\mathbb{R}^n$, in answering a long standing question of Payne-P\'olya-Weinberg \cite{ppw1}, Cheng and Yang  \cite{cy1} have obtained the following
universal inequality:
\be \sum_{i=1}^k(\laa_{k+1}-\laa_i)^2\leq
\fr{4(n+2)}{n^2}\sum_{i=1}^k(\laa_{k+1}-\laa_i)\laa_i.
\en
On the other hand, when $\Omega$  is a bounded domain in an $n$-dimensional unit sphere $\mathbb{S}^n$ and $l=2$, the following universal inequality has been proved in \cite{wx1}:
\be & &
2\sum_{i=1}^k(\laa_{k+1}-\laa_i)^2\\ \no &\leq&
\sum_{i=1}^k(\laa_{k+1}-\laa_i)^2\left(\delta \laa_i+\fr{\delta^2(\laa_i-(n-2))}{4(\delta\laa_i+n-2)}\right)  +\fr 1{\delta}\sum_{i=1}^k(\laa_{k+1}-\laa_i)\left(\laa_i+\fr{(n-2)^2}4\right),
\en
where $\delta$ is an arbitrary positive constant.

Recently, Cheng-Yang \cite{cy2} improved (1.2) to
\be \sum_{i=1}^k(\laa_{k+1}-\laa_i)^2\leq
\fr{4\left(n+\frac 43\right)}{n^2}\sum_{i=1}^k(\laa_{k+1}-\laa_i)\laa_i
\en
and conjectured that the following inequality is true:
\be \sum_{i=1}^k(\laa_{k+1}-\laa_i)^2\leq
\fr{4}{n}\sum_{i=1}^k(\laa_{k+1}-\laa_i)\laa_i.
\en
The inequality (1.3)  has been improved in \cite{jlwx1} and \cite{cy2}. For arbitrary $l$,  when $\Omega$ is a
bounded domain in a Euclidean space $\mathbb{R}^n$ or a unit sphere, Jost, Jost-Li,
Wang and Xia \cite{jlwx1} proved some universal inequalities which have been improved by Cheng, Qi, Wang and Xia \cite{cqwx}. For eigenvalues of the problem (1.1), in addition to considering the possible sharp inequalities, another interesting question is to
know for what kind of complete manifolds there exists universal upper bound on $\Lambda_{k+1}$ in terms of $\laa_1,\cdots, \laa_k$ independent of the domains.  Recently, universal inequalities for
eigenvalues of the (1.1) have been obtained  for bounded domains of some special Ricci flat manifolds in \cite{dwlx}.

In this paper, we prove the following result.
\vskip0.3cm
{\bf Theorem 1.1.} {\it Let $\mathbb{M}$ be a complete Riemannian manifold and let $\Omega $ be a bounded domain of the Riemannian product manifold $\mathbb{M}\times
\mathbb{R}$. Denote by  $\laa_i$  the ith eigenvalue of the problem (1.1). Then for any positive non-increasing monotone sequence $\{\delta_i\}_{i=1}^k$ we have
\be
& & \su \\ \no &\leq&\left(2l^2-\fr {11}3l+\fr 53\right)\sum_{i=1}^k \delta_i (\laa_{k+1}-\laa_i)^2 \laa_i^{(l-2)/(l-1)}+\sum_{i=1}^k\fr 1{\delta_i}(\laa_{k+1}-\laa_i)\laa_i^{1/(l-1)}.
\en
}
\vs
{\bf Corollary 1.2.} {\it Under the same conditions as in Theorem 1.1, we have}
\be
\su  \leq 4\left(2l^2-\fr {11}3l+\fr 53\right)\suu \laa_i.
\en

{\it Proof of Corollary 1.2.} Taking
\be\no
\delta_1=\cdots=\delta_k=\frac{\left(\sum_{i=1}^k(\laa_{k+1}-\laa_i)\laa_i^{1/(l-1)}\right)^{\fr 12}}{\left(\left(2l^2-\fr {11}3l+\fr 53\right)\sum_{i=1}^k(\laa_{k+1}-\laa_i)^2 \laa_i^{(l-2)/(l-1)}\right)^{\fr 12}}
\en
in (1.6) and using (Cf. \cite{jlwx2})
\be& &
\left(\su \laa_i^{(l-2)/(l-1)}\right)\left(\suu \laa_i^{1/(l-1)}\right)
\\ \no &\leq&
\left(\su \right)\left(\suu \laa_i\right),
\en
we get (1.7).

\section{A Proof of Theorem 1.1 }
\setcounter{equation}{0}
Before proving our result, let us recall some known facts we need.
Let  $\mathbb{M}$ be a complete manifold and let $\mathbb{N}=\mathbb{M}\times \mathbb{R}=\{(x, t)|x\in\mathbb{M}, t\in\mathbb{R}\}$ with the product metric of
$\mathbb{M}$ and $\mathbb{R}$. Let $\om(\subset\mathbb{N})$ be a bounded domain with smooth boundary. Since any complete Riemannian manifold can be isometrically embedded in some Euclidaen space, we can think of our $\mathbb{N}$ as a submanifold of some $\mathbb{R}^q$.
Let us denote by $\langle , \rangle$ the canonical metric on ${\mathbb R}^{q}$ as well as that induced on $\mathbb{N}$.
Denote by  $\Delta $ and $\nabla$    the Laplacian and the gradient operator of $\mathbb{N}$, respectively.
Let $u_i$ be the $i$-th orthonormal eigenfunction of the problem (1.1) corresponding to the eigenvalue $\laa_i$, $i=1, 2, \cdots,$ that is,
\be\left\{\begin{array}{lll}
 (-\Delta)^l u_i= -\laa_i\D u_i, \ \ {\rm in} \ \ \om, \\
   u_i= \fr{\pa u_i}{\pa \nu}=\cdots =\fr{\pa^{l-1} u_i}{\pa \nu^{l-1}}=0, \ \ {\rm on} \ \ \partial\om,\\
  (u_i, u_j)_D =\int_{\om}\langle \nabla u_i, \nabla u_j\rangle=\delta_{ij}, \ \forall\ i, j.
\end{array}\right.
\en
For functions $f$ and $g$ on $\om$, the {\it Dirichlet inner product $(f, g)_D$} of $f$ and $g$ is given by
\be \no (f, g)_D=\int_{\om}\langle\nabla f, \ \nabla g\rangle.
\en
The Dirichlet norm of a function $f$ is defined by
\be\no
||f||_D=\{(f, f)_D\}^{1/2}=\left(\int_{\om}|\nabla f|^2\right)^{1/2}.
\en

For each $k=1, \cdots, l$, let $\nabla^k$ denote the $k$-th covariant
derivative operator on $\mathbb{N}$, defined in the usual weak sense.
 For a function $f$ on $\om$, the squared norm of  $\nabla^k f$ is defined as (cf. \cite{he})
\be
\left|\nabla^kf\right|^2=\sum_{i_1,\cdots, i_k=1}^n\left(\nabla^kf(e_{i_1},\cdots, e_{i_k})\right)^2,
\en
where $e_1,\cdots, e_n$ are orthonormal vector fields locally defined on $\om$.
Define the Sobolev space $H_l^2(\om)$ by
$$H_l^2(\om)=\{ f:\ f, \ |\nabla f|,\cdots, \left|\nabla^l f\right|\in L^2(\om)\}.
$$
Then $H_l^2(\om)$ is a Hilbert space with respect to the inner product  $\langle\langle, \rangle\rangle$:
\be
\langle\langle f,  g \rangle\rangle=\int_{\om}\left(\sum_{k=0}^l \nabla^k f\cdot \na^k g\right),
\en
where
\be\no
\na^k f\cdot\na^k g = \sum_{i_1,\cdots, i_k=1}^n \nabla^kf(e_{i_1},\cdots, e_{i_k}) \nabla^kg(e_{i_1},\cdots, e_{i_k}).
\en
Consider the subspace $H_{l,D}^2(\om)$ of $H_l^2(\om)$ defined by
$$H_{l,D}^2(\om)=\left\{f\in H_l^2(\om): \ f|_{\pa \om}=\left. \fr{\pa f}{\pa \nu}\right|_{\pa \om}=\cdots\left. \fr{\pa^{l-1} f}{\pa \nu^{l-1}}\right|_{\pa \om}=0\right\}.
$$
The operator $(-\Delta)^l$ defines a self-adjoint operator acting on $H_{l,D}^2(\om)$
with discrete eigenvalues $0<\laa_1\leq\cdots\leq \laa_k\leq\cdots$ for the buckling problem (1.1) and the eigenfunctions $\{u_i\}_{i=1}^{\infty}$ defined in (2.1)
form a complete orthonormal basis for the Hilbert space $H_{l,D}^2(\om)$. If $\phi\in H_{l,D}^2(\om)$ satisfies $(\phi , u_j)_D=0, \ \forall  j=1, 2, \cdots, k$, then the Rayleigh-Ritz inequality tells us that
\be
\laa_{k+1}||\phi ||_D^2\leq \int_{\om} \phi(-\Delta)^l\phi.
\en
For vector-valued functions $F=(f_1, f_2, \cdots, f_{m}), \ G=(g_1, g_2, \cdots, g_{m}):
\om\rightarrow {\mathbb R}^{q}$, we define an inner product $(F, G)$ by
$$(F, G)\equiv \int_{\om} \langle F, G\rangle =\int_{\om} \sum_{\alpha =1}^{q} f_{\alpha}g_{\alpha}.$$
The norm of $F$ is given by
$$||F||=(F, F)^{1/2}=\left\{\int_{\om}\sum_{\alpha=1}^{q}f_{\alpha}^2\right\}^{1/2}.$$
Let ${\bf H}_1^2(\om)$ be the Hilbert space of vector-valued functions given by
$${\bf H}_1^2(\om)=\left\{ F=(f_1,\cdots, f_{q}): \om\rightarrow {\mathbb R}^{q};\  f_{\alpha}, \ |\nabla f_{\alpha}|\in L^2(\om), \ {\rm for} \ \alpha=1,\cdots, m\right\}
$$
with inner product $\langle, \rangle_1$:
$$\langle F, G\rangle_1=(F, G)+\int_{\om}\sum_{\alpha=1}^{q}\langle\nabla f_{\alpha}, \nabla g_{\alpha}\rangle.
$$
Observe that a vector field on $\om$ can be regarded as a vector-valued function from $\om$ to ${\mathbb R}^q$. Let ${\bf H}_{1, D}^2(\om)\subset  {\bf H}_{1}^2(\om)$ be a subspace of ${\bf H}_1^2(\om)$
spanned by the vector-valued functions $\{ \nabla u_i\}_{i=1}^{\infty}$, which form a complete orthonormal basis of ${\bf H}_{1, D}^2(\om)$. For any $f\in H_{l,D}^2(\om), $ we have $\nabla f\in {\bf H}_{1, D}^2(\om)$ and for any $X\in {\bf H}_{1, D}^2(\om)$, there exists a function $f\in H_{l,D}^2(\om)$ such that $X=\nabla f$.

Consider the function $g: \om\rightarrow \mathbb{R}$ given by $g(x, t)=t$. The vector fields $g\na u_i$ can be decomposed as
\be
g\na u_i= \na h_i +{\bf W}_i,
\en
where $h_i\in H_{l, D}^2(\om), \na h_i$ is the projection of $g\na h_i$ in ${\bf H}_{l, D}^2(\om )$, ${\bf W}_i \bot {\bf H}_{l, D}^2(\om )$ and
\be
\left\{\begin{array}{cl} {\bf W}_i\left|_{\pa \om}=0,\right.\\
{\rm div} {\bf W}_i =0.
\end{array}
\right.
\en
Since ${\mathbb N}$ is endowed with the product metric of ${\mathbb M}$ and $\mathbb R$, it is easy to see that
\be
|\na g|=1, \ \na^2g =0, \ {\rm Ric}(\na g, X)=0,
\en
for any vector fields $X\in {\mathfrak X}(\om)$, where ${\rm Ric}$ denotes the Ricci tensor of ${\mathbb N}$. Substituting (2.7) into the Bochner formula, we obtain for any smooth $f: \om\rightarrow\mathbb{R}$ that
\be
\D\lan\na f, \na g\ra &=&2\na^2f\cdot\na^2g+\lan\na f, \na(\D g)\ra +\lan\na g, \na(\D f)\ra +2{\rm Ric}(\na g, \na f)\\ \no
&=& \lan\na g, \na(\D f)\ra.
\en
In the proof of Theorem 1.1, we shall use (2.5)-(2.8) repeatedly.
Now we can prove the main result in this paper.
\vs
{\bf Proof of Theorem 1.1.}  For each $i=1,\cdots, k$, let us consider the functions $\phi_i: \om\ri  \mathbb R$ given by
\be
\phi_i=h_i-\sum_{j=1}^k b_{ij}u_j,
\en
where
$b_{ij}=\int_{\om} g\lan\na u_i, \na u_j\ra =b_{ji}$ and $h_i$ is given in (2.5).
Since
\be\no
\phi_i|_{\pa\om}=\left.\fr{\pa \phi_i}{\pa \nu}\right|_{\pa \om}=\cdots =\left.\fr{\pa^{l-1}\phi_i}{\pa \nu^{l-1}}\right|_{\pa \om}=0,
\int_{\om}\lan\na\phi_i, \na u_j\ra =0,\ j=1,\cdots,k,
\en
we have from the Rayleigh-Ritz inequality that
\be
\laa_{k+1}\int_{\om}|\na \phi_i|^2\leq \int_{\om}\phi_i(-\D)^l\phi_i, \ \forall i=1,\cdots,k.
\en
Using (2.5)-(2.9) and the same calculations as in \cite{dwlx} we get (Cf. (3.17)-(3.20) in \cite{dwlx})
\begin{eqnarray}& &
\ino \phi_i(-\D)^l\phi_i\\ \no &=&\ino
(-1)^{l}(2l^2-4l+3)\langle\na g,\na u_i\rangle\langle\na g,\na
(\D^{l-2}u_i)\rangle\\ \no&&+\ino
(-1)^{l}(-l+1)u_i\D^{l-1}u_i+\Li\left\{\ino g^2|\na u_i|^2-\ino
u_i^2\right\}-\sum_{j=1}^k\Lj b_{ij}^2,
\end{eqnarray}
\be
2\int_{\om}\left\lan g\na u_i, \na\lan \na g, \na u_i\ra \right\ra =-1,
\en
and
\be
c_{ij}\equiv \int_{\om}\left\lan\na\lan \na g, \na u_i\ra, \na u_j\right\ra = -c_{ji}.
\en
Substituting
\be
||g\na u_i||^2=||\na h_i||^2+||{\bf W}_i||^2, \ ||\na h_i||^2=||\na \phi_i||^2+\sum_{j=1}^k b^2_{ij}
\en
and (2.12) into (2.11), we get
\be& &
(\laa_{k+1}-\laa_i)||\na \phi_i||^2\\ \no &\leq&\ino
(-1)^{l}(2l^2-4l+3)\langle\na g,\na u_i\rangle\langle\na g,\na
(\D^{l-2}u_i)\rangle\\ \no&&+\ino
(-1)^{l}(-l+1)u_i\D^{l-1}u_i+\Li(||u_i||^2-||{\bf W}_i||^2)+\sum_{j=1}^k(\Li-\Lj) b_{ij}^2,
\en
where for a function $f: \om\ri {\mathbb R},\ ||f||^2=\int_{\om} f^2$ and for a vector field $X\in \mathfrak{X}(\om)$, $||X||^2=\int_{\om}|X|$, $|X|$ being the length of $X$.

It is easy to see from (2.5) and (2.9) that (Cf. (3.20) in \cite{dwlx})
\be
1+2\sum_{j=1}^k b_{ij}c_{ij}= -2\ino \left\lan\na\phi_i, \na \lan \na g, \na u_i\ra\right\ra.
\en
Thus, we have
\be& &
(\laa_{k+1}-\laa_i)^2\left(1+2\sum_{j=1}^k b_{ij}c_{ij}\right)
\\ \no &=&(\laa_{k+1}-\laa_i)^2\left(-2\ino \left\lan\na\phi_i, \na \lan \na g, \na u_i -\sum_{j=1}^k c_{ij}\na u_j\right\ra\right)\\ \no
&\leq&\delta_i (\laa_{k+1}-\laa_i)^3||\na \phi_i||^2+\fr 1{\delta_i}(\laa_{k+1}-\laa_i)\left(||\na\lan\na g, \na u_i\ra||^2-\sum_{j=1}^k c_{ij}^2\right).
\en
Summing on $i$ from $1$ to $k$ in (2.17),  using (2.15) and noticing $a_{ij}=a_{ji}, \ c_{ij}=-c_{ji}$, we get
\be& &
\su -2\sum_{i, j=1}^k (\laa_{k+1}-\laa_i)(\laa_i-\laa_j)b_{ij}c_{ij}\\ \no
&\leq&\sum_{i=1}^k \delta_i (\laa_{k+1}-\laa_i)^2\left(\ino
(-1)^{l}(2l^2-4l+3)\langle\na g,\na u_i\rangle\langle\na g,\na
(\D^{l-2}u_i)\rangle\right.\\ \no&&\left.+\ino
(-1)^{l}(-l+1)u_i\D^{l-1}u_i+\Li(||u_i||^2-||{\bf W}_i||^2)\right)
+\sum_{i=1}^k\fr 1{\delta_i}(\laa_{k+1}-\laa_i)||\na\lan\na g, \na u_i\ra||^2 \\ \no & &
+ \sum_{i, j=1}^k \delta_i (\laa_{k+1}-\laa_i)^2 (\laa_i-\laa_j)b_{ij}^2-\sum_{i, j=1}^k \fr 1{\delta_i}(\laa_{k+1}-\laa_i)c_{ij}^2.
\en
Since $\{\delta_i\}$ is a non-increasing monotone sequence, we have
\be& &
\sum_{i, j=1}^k \delta_i (\laa_{k+1}-\laa_i)^2 (\laa_i-\laa_j)b_{ij}^2\\ \no &=&-\sum_{i, j=1}^k \delta_i (\laa_{k+1}-\laa_i)(\laa_i-\laa_j)^2b_{ij}^2
+\fr 12\sum_{i, j=1}^k  (\laa_{k+1}-\laa_i)(\laa_{k+1}-\laa_j)(\delta_i-\delta_j) (\laa_i-\laa_j)b_{ij}^2\\ \no &\leq&
-\sum_{i, j=1}^k \delta_i (\laa_{k+1}-\laa_i)(\laa_i-\laa_j)^2b_{ij}^2.
\en
Substituting (2.19) into (2.18), we conclude
\be& &
\su \\ \no
&\leq&\sum_{i=1}^k \delta_i (\laa_{k+1}-\laa_i)^2\left(\ino
(-1)^{l}(2l^2-4l+3)\langle\na g,\na u_i\rangle\langle\na g,\na
(\D^{l-2}u_i)\rangle\right.\\ \no&&\left.+\ino
(-1)^{l}(-l+1)u_i\D^{l-1}u_i-\Li(||u_i||^2-||{\bf W}_i||^2)\right)
+\sum_{i=1}^k\fr 1{\delta_i}(\laa_{k+1}-\laa_i)||\na\lan\na g, \na u_i\ra||^2.
\en
Before we can finish the proof of Theorem 1.1, we shall need some lemmas.
\vs
{\bf Lemma 2.1.}  {\it We have}

 {\it i)} (Cf. (2.6) in \cite{jlwx1}). $0\leq \int_{\om} u_i (-\D )^k u_i \leq \laa_i^{(k-1)/(l-1)}, \ k=1,\cdots, l-1.$

 {\it ii)} (Cf. (3.26) in \cite{dwlx}). $||\na\lan \na g, \na u_i\ra ||^2\leq \laa_i^{1/(l-1)}$.

 {\it iii)} (Cf. (3.33) in \cite{dwlx}). $\int_{\om} \left\lan\na g, \na\left((-\D )^{l-2}u_i\right)\right\ra \lan \na g, \na u_i\ra \leq \laa_i^{(l-2)/(l-1)}$.
\vs
{\bf Lemma 2.2.} {\it We have}
\be
-\laa_i(||u_i||^2-||{\bf W}_i||^2) &\leq& -\fr{2(l+1)} 3 \int_{\om} \left\lan\na g, \na\left((-\D )^{l-2}u_i\right)\right\ra \lan \na g, \na u_i\ra
\\ \no & & +\fr 13\int_{\om} u_i(-\D)^{l-1}u_i.
\en
{\it Proof of Lemma 2.2.}
Firstly, we prove the following equality:
\be
-\laa_i(||u_i||^2-||{\bf W}_i||^2) &=& -\fr{l+1} 2 \int_{\om} \left\lan\na g, \na\left((-\D )^{l-2}u_i\right)\right\ra \lan \na g, \na u_i\ra
\\ \no & & -\fr 14\laa_i ||u_i||^2 +\fr 14\int_{\om} u_i(-\D)^{l-1}u_i.
\en
Observe that $\na(g u_i)=u_i\na g+ g\na u_i\in {\bf H}_{1, D}^2(\om )$. Set  $y_{i}= g u_i-h_{i}$; then
\be
u_i \na g=\na y_{i}-{\bf W}_{i}.
\en
and so
\be
||u_{i}||^2=||u_i \na g||^2=||{\bf W}_{i}||^2+||\na y_{i}||^2.
\en
We have
\be
- \ino y_i g \D u_i &=&  \ino \lan \na y_i, g\na u_i\ra + \ino y_i\lan\na g, \na u_i\ra
\\ \no &=&  \ino \lan \na y_i, g\na u_i\ra -\ino \lan\na y_i, u_i\na g\ra
\\ \no &=&  \ino \lan \na y_i, g\na u_i\ra -||\na y_i||^2\\ \no &= &
\ino \left\lan u_i\na g+{\bf W}_i, g\na u_i\right\ra -||\na y_i||^2\\ \no
&=& \fr 14 \ino \lan\na g^2, \na u_i^2\ra + ||{\bf W}_i||^2 -||\na y_i||^2\\ \no
&=& -\fr 14 \ino u_i^2 \D g^2  + ||{\bf W}_i||^2-||\na y_i||^2\\ \no &=& -\fr 12 ||u_i||^2 + ||{\bf W}_i||^2 -||\na y_i||^2
\\ \no &=& -2(||u_i||^2 - ||{\bf W}_i||^2)+\fr 12 ||u_i||^2,
\en
\be
\begin{aligned}
&\int_{\om}\D^{l-1}u_i\lan \na y_{i}, \na g\ra \\
=&-\int_{\om}y_{i}\lan \na (\D^{l-1}u_i), \na g\ra \\
=&-\int_{\om}y_{i}  \D\lan\na (\D^{l-2} u_i), \na g\ra \\
=&\int_{\om}\lan\na  y_{i}, \na\lan \na(\D^{l-2}u_i), \na g\ra\ra \\
=&\int_{\om}\lan u_i\na  g, \na\lan \na(\D^{l-2}u_i), \na g\ra\ra \\
=&-\int_{\om}\lan \na  g,  \na(\D^{l-2}u_i)\ra \lan\na u_i, \na g\ra
\end{aligned}
\en
and
\be
\begin{aligned}
&\int_{\om} g u_i\lan\na g, \na(\D^{l-1}u_i)\ra\\
=&\int_{\om} g u_i\D^{l-1}\lan\na g,\na u_i\ra\\
=&
\int_{\om} \D^{l-1}(g u_i)\lan\na g,\na u_i\ra\\
=&
-\int_{\om} \lan u_i\na g,\na(\D^{l-1}(g u_i))\ra\\
=&
-\int_{\om} \lan \na y_{i},\na(\D^{l-1}(g u_i))\ra\\
=&
\int_{\om}  y_{i}\D^{l}(g u_i)\\
=&
\int_{\om}  y_{i}(2l\lan\na(\D^{l-1}u_i), \na g\ra+g\D^lu_i)\\
=&
-2l\int_{\om}\D^{l-1}u_i\lan \na y_{i}, \na g\ra+\laa_i (-1)^{l-1}\int_{\om}y_{i}g\D u_i.
\end{aligned}
\en
It follows from (2.25)-(2.27) that
\be & &\int_{\om} g u_i\lan\na g, \na(\D^{l-1}u_i)\ra\\ \no &=&
2l\int_{\om}\lan \na  g,  \na(\D^{l-2}u_i)\ra \lan\na u_i, \na g\ra
+(-1)^{l-1}\laa_i\left(-\fr 12||u_i||^2+2(||u_i||^2-||{\bf W}_i||^2)\right).
\en
On the other hand, we also have
\be
\begin{aligned}
&\int_{\om} gu_i \lan\na g, \na(\D^{l-1}u_i) \ra \\
=&\int_{\om} gu_i \D^{l-1}\lan \na u_i, \na g\ra\\
=&\int_{\om}\D^{l-1}(gu_i) \lan \na u_i, \na g\ra\\
=&\int_{\om}\left(2(l-1)\lan \na (\D^{l-2}u_i), \na g\ra +g\D^{l-1} u_i\right) \lan \na u_i, \na g\ra
\end{aligned}
\en
and
\be
\int_{\om}gu_i \lan \na g, \na(\D^{l-1}u_i)\ra = -\int_{\om}\D^{l-1}u_i(u_i+ g\lan \na u_i, \na g\ra).
\en
Adding (2.29) and (2.30), we infer
\be & &
\int_{\om} gu_i \lan\na g, \na(\D^{l-1}u_i)\ra
 \\ \no &=&
 \int_M\left\{(l-1)(\lan \na(\D^{l-2}u_i), \na g\ra \lan \na u_i, \na g\ra-\fr 12u_i\D^{l-1}u_i \right\}.
\en
Combining (2.28) and (2.31), one gets (2.22).

Substituting
\be\no
-||u_i||^2\leq -(||u_i||^2-||{\bf W}_i|^2)
\en
into (2.22), we obtain (2.21). This completes the proof of Lemma 2.2.
\vs
Let us continue on the proof of Theorem 1.1. Substituting (2.21) into (2.20) and using Lemma 2.1, we have
\be& &
\su \\ \no
&\leq&\sum_{i=1}^k \delta_i (\laa_{k+1}-\laa_i)^2\left(\left(2l^2-4l+3-\fr{2(l+1)}3\right)\ino
\langle\na g,\na u_i\rangle\langle\na g,\na
((-\D)^{l-2}u_i)\rangle\right.\\ \no&&\left.+\ino
\left(l-1+\fr13\right)u_i(-\D)^{l-1}u_i\right)
+\sum_{i=1}^k\fr 1{\delta_i}(\laa_{k+1}-\laa_i)||\na\lan\na g, \na u_i\ra||^2\\ \no &\leq&
\sum_{i=1}^k \delta_i (\laa_{k+1}-\laa_i)^2 \left(2l^2-\frac{11}{3}l+\frac 53\right)\laa_i^{(l-2)/(l-1)}
\\ \no & & +\sum_{i=1}^k\fr 1{\delta_i}(\laa_{k+1}-\laa_i)\laa_i^{1/(l-1)}.
\en
This completes the proof of Theorem 1.1.
\section{Concluding Remarks}
By using (2.7), (2.8) and the same arguments as in the proof of Theorem 1.3 in \cite{dwlx}, one can prove the following result.
\vs
{\bf Theorem 3.1.} {\it Let $\mathbb{M}$ be a complete Riemannian manifold and let $\Omega $ be a bounded domain of the Riemannian product manifold $\mathbb{M}\times
\mathbb{R}$. Denote by  $\Gamma_i$ be the ith eigenvalue of the problem
\begin{eqnarray}\no
\left\{\begin{array}{l}(-\Delta)^lu =-\Gamma  u,~~\mbox{in}
~~\Omega,\\ ~~~~u=\frac{\partial u}{\partial
\nu}=\cdots=\frac{\partial^{l-1} u}{\partial
\nu^{n-1}}=0,~~\mbox{on}~~\partial \Omega.
\end{array}
\right.
\end{eqnarray}
Then we have
\be\no
\sum_{i=1}^k(\Gamma_{k+1}-\Gamma_i)^2\leq 2\left\{l(2l-1)\sum_{i=1}^k(\Gamma_{k+1}-\Gamma_i)\Gamma_i^{(l-1)/l}\right\}^{\frac 12}\left\{\sum_{i=1}^k(\Gamma_{k+1}-\Gamma_i)\Gamma_i^{
1/l}\right\}^{\frac 12}.
\en
}

\end{document}